# Dual Smarandache Curves and Smarandache Ruled Surfaces

**Tanju KAHRAMAN[a], Mehmet ÖNDER[a], H. Hüseyin UĞURLU[b]**

*[a] Celal Bayar University, Department of Mathematics, Faculty of Arts and Sciences, , Manisa, Turkey.*
*E-mails: tanju.kahraman@bayar.edu.tr, mehmet.onder@bayar.edu.tr*

*[b]Gazi University, Gazi Faculty of Education, Department of Secondary Education Science and Mathematics*
*Teaching, Mathematics Teaching Program, Ankara, Turkey. E-mail: hugurlu@gazi.edu.tr*

**Abstract**

In this paper, by considering dual geodesic trihedron (dual Darboux frame) we define dual Smarandache curves lying fully on dual unit sphere $\tilde{S}^2$ and corresponding to ruled surfaces. We obtain the relationships between the elements of curvature of dual spherical curve (ruled surface) $\tilde{\alpha}(s)$ and its dual Smarandache curve (Smarandache ruled surface) $\tilde{\alpha}_1(s)$ and we give an example for dual Smarandache curves of a dual spherical curve.

## 1. Introduction

In the Euclidean space $E^3$, an oriented line $L$ can be determined by a point $p \in L$ and a normalized direction vector $\vec{a}$ of $L$, i.e. $\|\vec{a}\| = 1$. The components of $L$ are obtained by the moment vector $\vec{a}^* = \vec{p} \times \vec{a}$ with respect to the origin in $E^3$. The two vectors $\vec{a}$ and $\vec{a}^*$ are not independent of one another; they satisfy the relationships $\langle \vec{a}, \vec{a} \rangle = 1$, $\langle \vec{a}, \vec{a}^* \rangle = 0$. The pair $(\vec{a}, \vec{a}^*)$ of the vectors $\vec{a}$ and $\vec{a}^*$, which satisfies those relationships, is called dual unit vector[2]. The most important properties of real vector analysis are valid for the dual vectors. Since each dual unit vector corresponds to a line of $E^3$, there is a one-to-one correspondence between the points of a dual unit sphere $\tilde{S}^2$ and the oriented lines of $E^3$. This correspondence is known as E. Study Mapping[2]. As a sequence of that, a differentiable curve lying fully on dual unit sphere in dual space $D^3$ represents a ruled surface which is a surface generated by moving of a line $L$ along a curve $\alpha(s)$ in $E^3$ and has the parametrization $\vec{r}(s,u) = \vec{\alpha}(s) + u\,\vec{l}(s)$, where $\vec{\alpha}(s)$ is called generating curve and $\vec{l}(s)$, the direction of the line $L$, is called ruling.

In the study of the fundamental theory and the characterizations of space curves, the special curves are very interesting and an important problem. The most mathematicians studied the special curves such as Mannheim curves and Bertrand curves. Recently, a new special curve which is called Smarandache curve is defined by Turgut and Yılmaz in Minkowski space-time[6]. Then Ali have studied Smarandache curves in the Euclidean 3-space $E^3$[1].

Moreover, Önder has studied the Bertrand offsets of ruled surface according to the dual geodesic trihedron(Darboux frame) and given the relationships between the dual and real curvatures of a ruled surface and its offset surface[5].

In this paper, we give Darboux approximation for dual Smarandache curves on dual unit sphere $\tilde{S}^2$. Firstly, we define the four types of dual Smarandache curves (Smarandache ruled surfaces) of a dual spherical curve(ruled surface). Then, we obtain the relationships between the dual curvatures of dual spherical curve $\tilde{\alpha}(s)$ and its dual Smarandache curves. Furthermore, we show that dual Smarandache $\tilde{e}\tilde{g}$ **-**curve of a dual curve is always its Bertrand offset. Finally, we give an example for Smarandache curves of an arbitrary curve on dual unit sphere $\tilde{S}^2$.

## 2. Dual Numbers and Dual Vectors

Let $D = IR \times IR = \{\bar{a} = (a, a^*) : a, a^* \in IR\}$ be the set of the pairs $(a, a^*)$. For $\bar{a} = (a, a^*)$, $\bar{b} = (b, b^*) \in D$ the following operations are defined on $D$:

Equality: $\bar{a} = \bar{b} \Leftrightarrow a = b,\ a^* = b^*$

Addition: $\bar{a} + \bar{b} = (a + b,\ a^* + b^*)$

Multiplication: $\bar{a}\bar{b} = (ab,\ ab^* + a^*b)$

The element $\varepsilon = (0,1) \in D$ satisfies the relationships

$$\varepsilon \neq 0, \quad \varepsilon^2 = 0, \quad \varepsilon 1 = 1\varepsilon = \varepsilon. \tag{1}$$

Let consider the element $\bar{a} \in D$ of the form $\bar{a} = (a, 0)$. Then the mapping $f : D \to IR,\ \ f(a,0) = a$ is a isomorphism. So, we can write $a = (a, 0)$. By the multiplication rule we have that

$$\begin{aligned} \bar{a} &= (a, a^*) \\ &= (a,0) + (0, a^*) \\ &= (a,0) + (0,1)(a^*, 0) \\ &= a + \varepsilon a^* \end{aligned} \tag{2}$$

Then $\bar{a} = a + \varepsilon a^*$ is called dual number and $\varepsilon$ is called dual unit. Thus the set of all dual numbers is given by

$$D = \{\bar{a} = a + \varepsilon a^* : a, a^* \in IR,\ \varepsilon^2 = 0\} \tag{3}$$

The set $D$ forms a commutative group under addition. The associative laws hold for multiplication. Dual numbers are distributive and form a ring over the real number field[2,4].

Dual function of dual number presents a mapping of a dual numbers space on itself. Properties of dual functions were thoroughly investigated by Dimentberg[3]. He derived the general expression for dual analytic (differentiable) function as follows

$$f(\bar{x}) = f(x + \varepsilon x^*) = f(x) + \varepsilon x^* f'(x), \tag{4}$$

where $f'(x)$ is derivative of $f(x)$ and $x, x^* \in IR$.

Let $D^3 = D \times D \times D$ be the set of all triples of dual numbers, i.e.,

$$D^3 = \{\tilde{a} = (\bar{a}_1, \bar{a}_2, \bar{a}_3) : \bar{a}_i \in D,\ i = 1,2,3\}, \tag{5}$$

Then the set $D^3$ is module together with addition and multiplication operations on the ring $D$ and called dual space. The elements of $D^3$ are called dual vectors. Similar to the dual numbers, a dual vector $\tilde{a}$ may be expressed in the form $\tilde{a} = \vec{a} + \varepsilon \vec{a}^* = (\vec{a}, \vec{a}^*)$, where $\vec{a}$ and $\vec{a}^*$ are the vectors of $IR^3$. Then for any vectors $\tilde{a} = \vec{a} + \varepsilon \vec{a}^*$ and $\tilde{b} = \vec{b} + \varepsilon \vec{b}^*$ of $D^3$, the scalar product and the cross product are defined by

$$\langle \tilde{a}, \tilde{b} \rangle = \langle \vec{a}, \vec{b} \rangle + \varepsilon\left(\langle \vec{a}, \vec{b}^* \rangle + \langle \vec{a}^*, \vec{b} \rangle\right), \tag{6}$$

and

$$\tilde{a} \times \tilde{b} = \vec{a} \times \vec{b} + \varepsilon\left(\vec{a} \times \vec{b}^* + \vec{a}^* \times \vec{b}\right), \tag{7}$$

respectively, where $\langle \vec{a}, \vec{b} \rangle$ and $\vec{a} \times \vec{b}$ are the inner product and the cross product of the vectors $\vec{a}$ and $\vec{a}^*$ in $IR^3$, respectively.

The norm of a dual vector $\tilde{a}$ is given by

$$\|\tilde{a}\| = \|\vec{a}\| + \varepsilon \frac{\langle \vec{a}, \vec{a}^* \rangle}{\|\vec{a}\|}, \quad (\vec{a} \neq 0). \tag{8}$$

A dual vector $\tilde{a}$ with norm $1+\varepsilon 0$ is called unit dual vector. The set of all dual unit vectors is given by

$$\tilde{S}^2 = \left\{ \tilde{a} = (a_1, a_2, a_3) \in D^3 : \langle \tilde{a}, \tilde{a} \rangle = 1 + \varepsilon 0 \right\}, \tag{9}$$

and called dual unit sphere[2,4].

E. Study used dual numbers and dual vectors in his research on the geometry of lines and kinematics. He devoted special attention to the representation of directed lines by dual unit vectors and defined the mapping that is known by his name:

***Theorem2.1.(E. Study Mapping):*** *There exists one-to-one correspondence between the vectors of dual unit sphere $\tilde{S}^2$ and the directed lines of space of lines $\mathbb{R}^3$* [2,4].

By the aid of this correspondence, the properties of the spatial motion of a line can be derived. Hence, the geometry of the ruled surface is represented by the geometry of dual curves lying fully on the dual unit sphere in $D^3$.

The angle $\bar{\theta} = \theta + \varepsilon\theta^*$ between two dual unit vectors $\tilde{a}, \tilde{b}$ is called *dual angle* and defined by

$$\langle \tilde{a}, \tilde{b} \rangle = \cos\bar{\theta} = \cos\theta - \varepsilon\theta^* \sin\theta.$$

By considering The E. Study Mapping, the geometric interpretation of dual angle is that $\theta$ is the real angle between lines $L_1$, $L_2$ corresponding to the dual unit vectors $\tilde{a}$, $\tilde{b}$ respectively, and $\theta^*$ is the shortest distance between those lines[2,4].

## 3. Dual Representation of Ruled Surfaces

In this section, we introduce dual representation of a ruled surface which is given by Veldkamp in [7] as follows:

Let $k$ be a dual curve represented by $\tilde{x} = \tilde{e}(u)$ or $\vec{x} + \varepsilon\vec{x}^* = \vec{e}(u) + \varepsilon\vec{e}^*(u)$. The real curve $x = e(u)$ on the real unit sphere is called *the (real) indicatrix of* $k$; we suppose throughout that it does not exist of a single point. We take as the parameter $u$ the arc-length $s$ on the real indicatrix and we denote differentiation with respect to $s$ by primes. Then $\tilde{x} = \tilde{e}(s)$ and $\langle \vec{e}', \vec{e}' \rangle = 1$. The vector $\vec{e}' = \vec{t}$ is the unit vector parallel to the tangent at the indicatrix. It is well known that given dual curve $k$ may be represented by

$$\tilde{x} = \tilde{e}(s) = \vec{e} + \varepsilon\vec{c} \times \vec{e}, \tag{10}$$

where

$$\langle \vec{e}, \vec{e} \rangle = 1, \ \langle \vec{e}', \vec{e}' \rangle = 1, \ \langle \vec{c}', \vec{e}' \rangle = 0.$$

We observe that $c$ is unambiguously determined by $k$. It follows from (10) that

$$\tilde{e}' = \vec{t} + \varepsilon\left(\vec{c} \times \vec{t} + \vec{c}' \times \vec{e}\right). \tag{11}$$

Hence by means of $|\tilde{x}| = |\vec{x}| + \varepsilon \frac{\langle \vec{x}, \vec{x}^* \rangle}{|\vec{x}|} \ \left(\vec{x} \neq 0\right)$:

$$\left|\tilde{e}'\right| = 1+\varepsilon \det\left(\vec{c}',\vec{e},\vec{t}\right) = 1+\varepsilon\Delta \tag{12}$$

where $\Delta = \det\left(\vec{c}',\vec{e},\vec{t}\right)$. Since $\vec{c}'$ as well as $\vec{e}$ is perpendicular to $\vec{t}$ we may write $\vec{c}'\times\vec{e} = \mu\vec{t}$; then we obtain $\Delta = \left\langle \vec{c}'\times\vec{e},\vec{t}\right\rangle = \mu$. Therefore $\vec{c}'\times\vec{e} = \Delta\vec{t}$ and we obtain in view of (11):

$$\tilde{e}' = \vec{t} + \varepsilon\left(\vec{c}\times\vec{t} + \Delta\vec{t}\right). \tag{13}$$

Let $\tilde{t}$ be dual unit vector with the same sense as $\tilde{e}'$; then we find as a consequence of (12): $\tilde{e}' = (1+\varepsilon\Delta)\tilde{t}$. This leads in view of (13) to:

$$\tilde{t} = \vec{t} + \varepsilon\vec{c}\times\vec{t}. \tag{14}$$

Guided by elementary differential geometry of real curves we introduce the *dual arc-length* $\bar{s}$ of the dual curve $k$ by

$$\bar{s} = \int_0^s \left|\tilde{e}'(\sigma)\right| d\sigma = \int_0^s (1+\varepsilon\Delta)d\sigma = s + \varepsilon\int_0^s \Delta d\sigma.$$

Then $\bar{s}' = 1+\varepsilon\Delta$. We define furthermore: $\dfrac{d\tilde{e}}{d\bar{s}} = \dfrac{\tilde{e}'}{\bar{s}'}$; hence $\dfrac{d\tilde{e}}{d\bar{s}} = \tilde{e}' + \varepsilon\Delta$ and therefore

$$\frac{d\tilde{e}}{d\bar{s}} = \tilde{t}. \tag{15}$$

Introducing the dual unit vector $\tilde{e}\times\tilde{t} = \tilde{g} = \vec{g} + \varepsilon\vec{g}^*$ we observe $\vec{e}\times\vec{t} = \vec{g}$; hence by means of (10) and (14):

$$\tilde{g} = \vec{g} + \varepsilon\vec{c}\times\vec{g}. \tag{16}$$

Then the dual frame $\left\{\tilde{e},\tilde{t},\tilde{g}\right\}$ is called dual geodesic trihedron( or dual Darboux frame) of the ruled surface corresponding to dual curve $\tilde{e}$. Thus, the derivative formulae of this frame are given as follows,

$$\frac{d\tilde{e}}{d\bar{s}} = \tilde{t},\ \frac{d\tilde{t}}{d\bar{s}} = \bar{\gamma}\,\tilde{g} - \tilde{e},\ \frac{d\tilde{g}}{d\bar{s}} = -\bar{\gamma}\,\tilde{t} \tag{17}$$

where $\bar{\gamma}$ is called dual spherical curvature and given by

$$\bar{\gamma} = \gamma + \varepsilon\left(\delta - \gamma\Delta\right); \tag{18}$$

and $\delta = \left\langle \vec{c}',\vec{e}\right\rangle$, $\gamma = -\left\langle \vec{g}',\vec{t}\right\rangle$. From (17) introducing the *dual Darboux vector* $\tilde{d} = \bar{\gamma}\,\tilde{e} + \tilde{g}$ we have

$$\frac{d\tilde{e}}{d\bar{s}} = \tilde{d}\times\tilde{e},\ \frac{d\tilde{t}}{d\bar{s}} = \tilde{d}\times\tilde{t},\ \frac{d\tilde{g}}{d\bar{s}} = \tilde{d}\times\tilde{g}. \tag{19}$$

(See [8]).

Analogous to common differential geometry the dual radius of curvature $\bar{R}$ of the dual curve $\tilde{x} = \tilde{e}(s)$ is given by

$$\bar{R} = \frac{\left|\dfrac{d\tilde{e}}{d\bar{s}}\right|^3}{\left|\dfrac{d\tilde{e}}{d\bar{s}}\times\dfrac{d^2\tilde{e}}{d\bar{s}^2}\right|}.$$

Then from (15) and (17),

$$\bar{R} = \left(1+\bar{\gamma}^2\right)^{-1/2}. \tag{20}$$

The unit vector $\tilde{d}_0$ with the same sense as the Darboux vector $\tilde{d} = \bar{\gamma}\tilde{e} + \tilde{g}$ is given by

$$\tilde{d}_0 = \frac{\bar{\gamma}}{\sqrt{1+\bar{\gamma}^2}}\tilde{e} + \frac{1}{\sqrt{1+\bar{\gamma}^2}}\tilde{g}. \tag{21}$$

The dual angle $\bar{\rho}$ between $\tilde{d}_0$ and $\tilde{e}$ satisfies therefore:

$$\cos\bar{\rho} = \frac{\bar{\gamma}}{\sqrt{1+\bar{\gamma}^2}} \quad , \quad \sin\bar{\rho} = \frac{1}{\sqrt{1+\bar{\gamma}^2}} .$$

Hence:

$$\bar{R} = \sin\bar{\rho} \ , \quad \bar{\gamma} = \cot\bar{\rho} . \tag{22}$$

The point $\bar{M}$ on the dual unit sphere indicated by $\tilde{d}_0$ is called the *dual spherical centre of curvature* of $k$ at the point $Q$ given by the parameter value $\bar{s}$, whereas $\bar{\rho}$ is the *dual spherical radius of curvature*[7].

## 5. Dual Smarandache curves and Smarandache Ruled Surfaces

From E. Study Mapping, it is well-known that dual curves lying on dual unit sphere correspond to ruled surfaces of the line space $IR^3$. Thus, by defining the dual smarandache curves lying fully on dual unit sphere, we also define the smarandache ruled surfaces. Then, the differential geometry of smarandache ruled surfaces can be investigated by considering the corresponding dual smarandache curves on dual unit sphere.

In this section, we first define the four different types of the dual smarandache curves on dual unit sphere. Then by the aid of dual geodesic trihedron(Dual Darboux frame), we give the characterizations of these dual curves(or ruled surfaces).

### 5.1. Dual Smarandache $\tilde{e}\tilde{t}$ -curve of a unit dual spherical curve(ruled surface)

In this section, we define the first type of dual Smarandache curves as dual Smarandache $\tilde{e}\tilde{t}$ -curve. Then, we give the relationships between the dual curve and its dual Smarandache $\tilde{e}\tilde{t}$ -curve. Using the found results and relationships we study the developable of the corresponding ruled surface and its Smaranadache ruled surface.

**Definition 5.1.** Let $\tilde{\alpha} = \tilde{\alpha}(\bar{s})$ be a unit speed regular dual curve lying fully on dual unit sphere $\tilde{S}^2$ and $\{\tilde{e}, \tilde{t}, \tilde{g}\}$ be its moving dual Darboux frame. The dual curve $\tilde{\alpha}_1$ defined by

$$\tilde{\alpha}_1 = \frac{1}{\sqrt{2}}\left(\tilde{e} + \tilde{t}\right) \tag{23}$$

is called the dual Smarandache $\tilde{e}\tilde{t}$ -curve of $\tilde{\alpha}$ and fully lies on $\tilde{S}^2$. Then the ruled surface corresponding to $\tilde{\alpha}_1$ is called the Smarandache $\vec{e}\vec{t}$ -ruled surface of the surface corresponding to dual curve $\tilde{\alpha}$.

Now we can give the relationships between $\tilde{\alpha}$ and its dual Smarandache $\tilde{e}\tilde{t}$ -curve $\tilde{\alpha}_1$ as follows.

***Theorem 5.1.*** *Let $\tilde{\alpha} = \tilde{\alpha}(\bar{s})$ be a unit speed regular dual curve lying on dual unit sphere $\tilde{S}^2$. Then the relationships between the dual Darboux frames of $\tilde{\alpha}$ and its dual Smarandache $\tilde{e}\tilde{t}$ -curve $\tilde{\alpha}_1$ are given by*

$$\begin{pmatrix} \tilde{e}_1 \\ \tilde{t}_1 \\ \tilde{g}_1 \end{pmatrix} = \begin{pmatrix} \frac{1}{\sqrt{2}} & \frac{1}{\sqrt{2}} & 0 \\ -\frac{1}{\sqrt{2+\overline{\gamma}^2}} & \frac{1}{\sqrt{2+\overline{\gamma}^2}} & \frac{\overline{\gamma}}{\sqrt{2+\overline{\gamma}^2}} \\ \frac{\overline{\gamma}}{\sqrt{4+2\overline{\gamma}^2}} & -\frac{\overline{\gamma}}{\sqrt{4+2\overline{\gamma}^2}} & \frac{\sqrt{2}}{\sqrt{2+\overline{\gamma}^2}} \end{pmatrix} \begin{pmatrix} \tilde{e} \\ \tilde{t} \\ \tilde{g} \end{pmatrix} \tag{24}$$

*where $\overline{\gamma}$ is as given in (18).*

**Proof.** Let us investigate the dual Darboux frame fields of dual Smarandache $\tilde{e}\tilde{t}$-curve according to $\tilde{\alpha} = \tilde{\alpha}(\overline{s})$. Since $\tilde{\alpha}_1 = \tilde{e}_1$, we have

$$\tilde{e}_1 = \frac{1}{\sqrt{2}}\left(\tilde{e} + \tilde{t}\right) \tag{25}$$

Differentiating (25) with respect to $\overline{s}$, we get

$$\frac{d\tilde{e}_1}{d\overline{s}} = \frac{d\tilde{e}_1}{d\overline{s}_1} \cdot \frac{d\overline{s}_1}{d\overline{s}} = \tilde{t}_1 \cdot \frac{d\overline{s}_1}{d\overline{s}} = \frac{1}{\sqrt{2}}\left(-\tilde{e} + \tilde{t} + \overline{\gamma}\tilde{g}\right)$$

and hence

$$\tilde{t}_1 = \frac{\left(-\tilde{e} + \tilde{t} + \overline{\gamma}\tilde{g}\right)}{\sqrt{2+\overline{\gamma}^2}} \tag{26}$$

where

$$\frac{d\overline{s}_1}{d\overline{s}} = \sqrt{\frac{2+\overline{\gamma}^2}{2}}.$$

Thus, since $\tilde{g}_1 = \tilde{e}_1 \times \tilde{t}_1$, we have

$$\tilde{g}_1 = \frac{\overline{\gamma}}{\sqrt{4+2\overline{\gamma}^2}}\,\tilde{e} - \frac{\overline{\gamma}}{\sqrt{4+2\overline{\gamma}^2}}\,\tilde{t} + \frac{\sqrt{2}}{\sqrt{2+\overline{\gamma}^2}}\,\tilde{g}\,. \tag{27}$$

From (25)-(27) we have (24).

If we represent the dual Darboux frames of $\tilde{\alpha}$ and $\tilde{\alpha}_1$ by the dual matrixes $\tilde{E}$ and $\tilde{E}_1$, respectively, then (24) can be written as follows

$$\tilde{E} = \tilde{A}\tilde{E}_1$$

where

$$\tilde{A} = \begin{pmatrix} \frac{1}{\sqrt{2}} & \frac{1}{\sqrt{2}} & 0 \\ -\frac{1}{\sqrt{2+\overline{\gamma}^2}} & \frac{1}{\sqrt{2+\overline{\gamma}^2}} & \frac{\overline{\gamma}}{\sqrt{2+\overline{\gamma}^2}} \\ \frac{\overline{\gamma}}{\sqrt{4+2\overline{\gamma}^2}} & -\frac{\overline{\gamma}}{\sqrt{4+2\overline{\gamma}^2}} & \frac{\sqrt{2}}{\sqrt{2+\overline{\gamma}^2}} \end{pmatrix}. \tag{28}$$

It is easily seen that $\det(\tilde{A}) = 1$ and $\tilde{A}\tilde{A}^T = \tilde{A}^T\tilde{A} = I$ where $I$ is the $3\times 3$ unitary matrix. It means that $\tilde{A}$ is a dual orthogonal matrix. Then we can give the following corollary.

***Corollary 5.1.*** *The relationship between Darboux frames of the dual curves(ruled surfaces) $\tilde{\alpha}$ and $\tilde{\alpha}_1$ is given by a dual orthogonal matrix defined in (28).*

***Theorem 5.2.*** *The relationship between the dual Darboux formulae of dual Smarandache $\tilde{e}\tilde{t}$ -curve $\tilde{\alpha}_1$ and dual Darboux frame of $\tilde{\alpha}$ is as follows*

$$\begin{pmatrix} \dfrac{d\tilde{e}_1}{d\bar{s}_1} \\ \dfrac{d\tilde{t}_1}{d\bar{s}_1} \\ \dfrac{d\tilde{g}_1}{d\bar{s}_1} \end{pmatrix} = \begin{pmatrix} \dfrac{-1}{\sqrt{2+\bar{\gamma}^2}} & \dfrac{1}{\sqrt{2+\bar{\gamma}^2}} & \dfrac{\bar{\gamma}}{\sqrt{2+\bar{\gamma}^2}} \\ \dfrac{\sqrt{2}\bar{\gamma}\bar{\gamma}'-\sqrt{2}\bar{\gamma}^2-2\sqrt{2}}{\left(2+\bar{\gamma}^2\right)^2} & \dfrac{-\sqrt{2}\bar{\gamma}\bar{\gamma}'-\sqrt{2}(1+\bar{\gamma}^2)(2+\bar{\gamma}^2)}{\left(2+\bar{\gamma}^2\right)^2} & \dfrac{\sqrt{2}(\bar{\gamma}'+\bar{\gamma})(2+\bar{\gamma}^2)-\sqrt{2}\bar{\gamma}^2\bar{\gamma}'}{\left(2+\bar{\gamma}^2\right)^2} \\ \dfrac{\bar{\gamma}^3+2\bar{\gamma}'+2\bar{\gamma}}{\left(2+\bar{\gamma}^2\right)^2} & \dfrac{-\bar{\gamma}^3-2\bar{\gamma}'-2\bar{\gamma}}{\left(2+\bar{\gamma}^2\right)^2} & \dfrac{-\bar{\gamma}^4-2\bar{\gamma}^2-2\bar{\gamma}\bar{\gamma}'}{\left(2+\bar{\gamma}^2\right)^2} \end{pmatrix} \begin{pmatrix} \tilde{e} \\ \tilde{t} \\ \tilde{g} \end{pmatrix} \tag{29}$$

**Proof.** Differentiating (25), (26) and (27) with respect to $\bar{s}$, we have the desired equation (29).

***Theorem 5.3.*** *Let $\tilde{\alpha}=\tilde{\alpha}(\bar{s})$ be a unit speed regular curve on dual unit sphere. Then the relationship between the dual curvatures of $\tilde{\alpha}$ and its dual Smarandache $\tilde{e}\tilde{t}$ -curve $\tilde{\alpha}_1$ is given by*

$$\bar{\gamma}_1=\frac{\bar{\gamma}^3+2\bar{\gamma}'+2\bar{\gamma}}{\left(2+\bar{\gamma}^2\right)^{\frac{3}{2}}}. \tag{30}$$

**Proof.** Since $\dfrac{d\tilde{g}_1}{d\bar{s}_1}=-\bar{\gamma}_1\tilde{t}_1$, from (26) and (28), we get dual curvature of the curve $\tilde{\alpha}_1(\bar{s}_1)$ as follows

$$\bar{\gamma}_1=\frac{\bar{\gamma}^3+2\bar{\gamma}'+2\bar{\gamma}}{\left(2+\bar{\gamma}^2\right)^{\frac{3}{2}}}.$$

***Corollary 5.2.*** *If the dual curvature $\bar{\gamma}$ of $\tilde{\alpha}$ is zero, then the dual curvature $\bar{\gamma}_1$ of dual Smarandache $\tilde{e}\tilde{t}$ -curve $\tilde{\alpha}_1$ is zero.*

***Corollary 5.3.*** *The Darboux instantaneous vector of dual Smarandache $\tilde{e}\tilde{t}$ -curve is given by*

$$\tilde{d}_1=\frac{1}{\sqrt{2}\left(2+\bar{\gamma}^2\right)^{\frac{3}{2}}}\left[(2\bar{\gamma}^3+2\bar{\gamma}'+4\bar{\gamma})\,\tilde{e}+2\bar{\gamma}'\,\tilde{t}+(2\bar{\gamma}^2+4)\,\tilde{g}\right]. \tag{31}$$

**Proof:** It is known that the dual Darboux instantaneous vector of dual Smarandache $\tilde{e}\tilde{t}$ -curve is $\tilde{d}_1=\bar{\gamma}_1\tilde{e}_1+\tilde{g}_1$. Then, from (24) and (30) we have (31).

***Theorem 5.4.*** *Let $\tilde{\alpha}(\bar{s})=\tilde{\alpha}$ be a unit speed regular dual curve on dual unit sphere and $\tilde{\alpha}_1$ be its dual Smarandache $\tilde{e}\tilde{t}$ -curve. If the ruled surface corresponding to dual curve $\tilde{\alpha}$ is developable then the ruled surface corresponding to dual curve $\tilde{\alpha}_1$ is also developable if and only if*

$$\delta_1 = \frac{3\delta\gamma^2 + 2\delta + 2\delta'}{\left(2+\gamma^2\right)^{\frac{3}{2}}} - \frac{3\delta\gamma(\gamma^3 + 2\gamma' + 2\gamma)}{\left(2+\gamma^2\right)^{\frac{5}{2}}}$$

**Proof.** From (18) we have

$$\overline{\gamma} = \gamma + \varepsilon\left(\delta - \gamma\Delta\right)$$

$$\overline{\gamma}_1 = \gamma_1 + \varepsilon\left(\delta_1 - \gamma_1\Delta_1\right)$$

Then substituting these equalities into to equation (30), we have

$$\delta_1 - \gamma_1\Delta_1 = \frac{(3\gamma^2+2)(\delta - \gamma\Delta) + 2\delta' - 2\gamma'\Delta - 2\gamma\Delta'}{\left(2+\gamma^2\right)^{\frac{3}{2}}} - \frac{3\gamma(\delta - \gamma\Delta)(\gamma^3 + 2\gamma' + 2\gamma)}{\left(2+\gamma^2\right)^{\frac{5}{2}}}.$$

Since the ruled surface corresponding to dual curve $\tilde{\alpha}$ is developable, $\Delta = 0$. Hence,

$$\Delta_1 = \frac{\delta_1}{\gamma_1} - \frac{(3\gamma^2+2)\delta + 2\delta'}{\gamma_1\left(2+\gamma^2\right)^{\frac{3}{2}}} + \frac{3\gamma\delta(\gamma^3 + 2\gamma' + 2\gamma)}{\gamma_1\left(2+\gamma^2\right)^{\frac{5}{2}}}.$$

Thus, the ruled surface corresponding to dual curve $\tilde{\alpha}_1$ is developable if and only if

$$\delta_1 = \frac{3\delta\gamma^2 + 2\delta + 2\delta'}{\left(2+\gamma^2\right)^{\frac{3}{2}}} - \frac{3\delta\gamma(\gamma^3 + 2\gamma' + 2\gamma)}{\left(2+\gamma^2\right)^{\frac{5}{2}}}.$$

***Theorem 5.5.*** *The relationship between the radius of dual curvature of dual Smarandache $\tilde{e}\tilde{t}$-curve $\tilde{\alpha}_1$ and the dual curvature of $\tilde{\alpha}$ is given by*

$$\overline{R}_1 = \frac{\left(2+\overline{\gamma}^2\right)^{\frac{3}{2}}}{\sqrt{2\overline{\gamma}^6 + 14\overline{\gamma}^4 + 12\overline{\gamma}^2 + 4\overline{\gamma}^3\overline{\gamma}' + 8\overline{\gamma}\overline{\gamma}' + 4\overline{\gamma}'^2}}. \quad (32)$$

**Proof.** From (20), $\overline{R}_1 = \dfrac{1}{\left(1+\overline{\gamma}_1^2\right)^{2}}$. Then from (30), the radius of dual curvature is

$$\overline{R}_1 = \frac{1}{\sqrt{1 + \dfrac{\left(\overline{\gamma}^3 + 2\overline{\gamma}' + 2\overline{\gamma}\right)^2}{\left(2+\overline{\gamma}^2\right)^3}}} = \frac{\left(2+\overline{\gamma}^2\right)^{\frac{3}{2}}}{\sqrt{2\overline{\gamma}^6 + 14\overline{\gamma}^4 + 12\overline{\gamma}^2 + 4\overline{\gamma}^3\overline{\gamma}' + 8\overline{\gamma}\overline{\gamma}' + 4\overline{\gamma}'^2}}.$$

***Theorem 5.6.*** *The relationship between the radius of dual spherical curvature of dual Smarandache $\tilde{e}\tilde{t}$-curve $\tilde{\alpha}_1$ and the elements of dual curvature of $\tilde{\alpha}$ is,*

$$\overline{\rho}_1 = \arcsin\left(\frac{\left(2+\overline{\gamma}^2\right)^{\frac{3}{2}}}{\sqrt{2\overline{\gamma}^6 + 14\overline{\gamma}^4 + 12\overline{\gamma}^2 + 4\overline{\gamma}^3\overline{\gamma}' + 8\overline{\gamma}\overline{\gamma}' + 4\overline{\gamma}'^2}}\right) \quad (33)$$

**Proof.** Let $\overline{\rho}_1$ be the radius of dual spherical curvature and $\overline{R}_1$ be the radius of dual curvature of $\tilde{\alpha}_1$. From equation (22) we have

$$\sin\overline{\rho}_1 = \overline{R}_1$$

Thus, we get radius of dual spherical curvature

$$\bar{\rho}_1 = \arcsin\left(\frac{\left(2+\bar{\gamma}^2\right)^{\frac{3}{2}}}{\sqrt{2\bar{\gamma}^6+14\bar{\gamma}^4+12\bar{\gamma}^2+4\bar{\gamma}^3\bar{\gamma}'+8\bar{\gamma}\bar{\gamma}'+4\bar{\gamma}'^2}}\right).$$

In the following sections we define dual Smarandache $\tilde{e}\tilde{g}$, $\tilde{t}\tilde{g}$ and $\tilde{e}\tilde{t}\tilde{g}$ curves. The proofs of the theorems and corollaries of these sections can be given by using the similar way used in previous section.

**5.2. Dual Smarandache $\tilde{e}\tilde{g}$-curve of a unit dual spherical curve(ruled surface)**

In this section, we define the second type of dual Smarandache curves as dual Smarandache $\tilde{e}\tilde{g}$-curve. Then, we give the relationships between the dual curve and its dual Smarandache $\tilde{e}\tilde{g}$-curve. Using obtained results and relationships we study the developable of the corresponding ruled surface and its Smarandache ruled surface.

**Definition 5.2.** Let $\tilde{\alpha}(\bar{s})=\tilde{\alpha}$ be a unit speed regular dual curve lying fully on dual unit sphere and $\{\tilde{e},\tilde{t},\tilde{g}\}$ be its moving Darboux frame. The dual curve $\tilde{\alpha}_2$ defined by

$$\tilde{\alpha}_2 = \frac{1}{\sqrt{2}}(\tilde{e}+\tilde{g}). \tag{34}$$

is called the dual Smarandache $\tilde{e}\tilde{g}$-curve of $\tilde{\alpha}$ and fully lies on $\tilde{S}^2$. Then the ruled surface corresponding to $\tilde{\alpha}_2$ is called the Smarandache $\tilde{e}\tilde{g}$-ruled surface of the surface corresponding to dual curve $\tilde{\alpha}$.

Now we can give the relationships between $\tilde{\alpha}$ and its dual smarandache $\tilde{e}\tilde{g}$-curve $\tilde{\alpha}_2$ as follows.

***Theorem 5.7.*** *Let $\tilde{\alpha}(\bar{s})=\tilde{\alpha}$ be a unit speed regular dual curve lying on dual unit sphere $\tilde{S}^2$. Then the relationships between the dual Darboux frames of $\tilde{\alpha}$ and its dual Smarandache $\tilde{e}\tilde{g}$-curve $\tilde{\alpha}_2$ are given by*

$$\begin{pmatrix}\tilde{e}_2\\ \tilde{t}_2\\ \tilde{g}_2\end{pmatrix} = \begin{pmatrix}\frac{1}{\sqrt{2}} & 0 & \frac{1}{\sqrt{2}}\\ 0 & 1 & 0\\ \frac{-1}{\sqrt{2}} & 0 & \frac{1}{\sqrt{2}}\end{pmatrix}\begin{pmatrix}\tilde{e}\\ \tilde{t}\\ \tilde{g}\end{pmatrix}. \tag{35}$$

From (35) we have $\tilde{t}_2=\tilde{t}$, i.e, $\tilde{\alpha}_2$ is a Bertrand offset of $\tilde{\alpha}$ [5].

In [5], Önder has given the relationship between the geodesic frames of Bertrand surface offsets as follows

$$\begin{pmatrix}\tilde{e}_2\\ \tilde{t}_2\\ \tilde{g}_2\end{pmatrix} = \begin{pmatrix}\cos\bar{\theta} & 0 & -\sin\bar{\theta}\\ 0 & 1 & 0\\ \sin\bar{\theta} & 0 & \cos\bar{\theta}\end{pmatrix}\begin{pmatrix}\tilde{e}\\ \tilde{t}\\ \tilde{g}\end{pmatrix}$$

where $\bar{\theta}=\theta+\varepsilon\theta^*$, $(0\le\theta\le\pi,\ \theta^*\in\mathbb{R})$ is the dual angle between the generators $\tilde{e}$ and $\tilde{e}_2$ of Bertrand ruled surface $\varphi_e$ and $\varphi_{e_2}$. The angle $\theta$ is called the offset angle and $\theta^*$ is called the offset distance[5]. Then from (35) we have that offset angle is $\theta=\pi/4$ and offset distance is $\theta^*=0$. Then we have the following corollary.

***Corollary 5.4.*** *The dual Smarandache* $\tilde{e}\tilde{g}$ *-curve of a dual curve* $\tilde{\alpha}$ *is always its Bertrand offset with dual offset angle* $\overline{\theta} = \pi / 4 + \varepsilon 0$.

***Theorem 5.8.*** *Let* $\tilde{\alpha}(\overline{s}) = \tilde{\alpha}$ *be a unit speed regular dual curve on dual unit sphere* $\tilde{S}^2$. *Then according to dual Darboux frame of* $\tilde{\alpha}$, *the dual Darboux formulae of dual Smarandache* $\tilde{e}\tilde{g}$ *-curve* $\tilde{\alpha}_2$ *are as follows*

$$\begin{pmatrix} \dfrac{d\tilde{e}_2}{ds_2} \\ \dfrac{d\tilde{t}_2}{ds_2} \\ \dfrac{d\tilde{g}_2}{ds_2} \end{pmatrix} = \begin{pmatrix} 0 & 1 & 0 \\ \dfrac{-\sqrt{2}}{1-\overline{\gamma}} & 0 & \dfrac{\sqrt{2}\,\overline{\gamma}}{1-\overline{\gamma}} \\ 0 & -\dfrac{1+\overline{\gamma}}{1-\overline{\gamma}} & 0 \end{pmatrix} \begin{pmatrix} \tilde{e} \\ \tilde{t} \\ \tilde{g} \end{pmatrix} \tag{36}$$

***Theorem 5.9.*** *Let* $\tilde{\alpha}(\overline{s}) = \tilde{\alpha}$ *be a unit speed regular curve on dual unit sphere. Then the relationship between the dual curvatures of* $\tilde{\alpha}$ *and its dual Smarandache* $\tilde{e}\tilde{g}$ *-curve* $\tilde{\alpha}_2$ *is given by*

$$\overline{\gamma}_2 = \frac{1+\overline{\gamma}}{1-\overline{\gamma}} \ .$$

***Corollary 5.5.*** *The dual curvature* $\overline{\gamma}$ *of* $\tilde{\alpha}$ *is zero if and only if the dual curvature* $\overline{\gamma}_2$ *of dual Smarandache* $\tilde{e}\tilde{t}$ *-curve* $\tilde{\alpha}_1$ *is* 1.

***Corollary 5.6.*** *The Darboux instantaneous vector of dual Smarandache* $\tilde{e}\tilde{g}$ *-curve is given by*

$$\tilde{d}_2 = \frac{\sqrt{2}\overline{\gamma}}{1-\overline{\gamma}}\,\tilde{e} + \frac{\sqrt{2}}{1-\overline{\gamma}}\,\tilde{g} \ .$$

***Theorem 5.10.*** *Let* $\tilde{\alpha}(\overline{s}) = \tilde{\alpha}$ *be a unit speed regular curve on dual unit sphere and* $\tilde{\alpha}_2$ *be the dual Smarandache* $\tilde{e}\tilde{g}$ *-curve of* $\tilde{\alpha}$. *If the ruled surface corresponding to the dual curve* $\tilde{\alpha}$ *is developable then the ruled surface corresponding to dual curve* $\tilde{\alpha}_2$ *is developable if and only if*

$$\left(1-\gamma^2\right)\delta_2 - 2\delta = 0 .$$

***Theorem 5.11.*** *Let* $\tilde{\alpha}(\overline{s}) = \tilde{\alpha}$ *be a unit speed regular curve on dual unit sphere. Then the relationship between the radius of dual curvature of dual Smarandache* $\tilde{e}\tilde{g}$ *-curve* $\tilde{\alpha}_2$ *and the dual curvature of* $\tilde{\alpha}(\overline{s}) = \tilde{\alpha}$ *is,*

$$\overline{R}_2 = \frac{1-\overline{\gamma}}{\sqrt{2\overline{\gamma}^2 + 2}} .$$

***Theorem 5.12.*** *Let* $\tilde{\alpha}(\overline{s}) = \tilde{\alpha}$ *be a unit speed regular curve on dual unit sphere. Then the relationship between the radius of dual spherical curvature of dual Smarandache* $\tilde{e}\tilde{g}$ *-curve* $\tilde{\alpha}_2$ *and the elements of dual curvature of* $\tilde{\alpha}$ *is,*

$$\bar{\rho}_2 = \arcsin\left(\frac{1-\gamma}{\sqrt{2+2\gamma^2}}\right) - \varepsilon \frac{\gamma^* + \gamma\gamma^*}{\sqrt{2}\left(1+\gamma^2\right)^{\frac{3}{2}} \cos\left(\arcsin\left(\frac{1-\gamma}{\sqrt{2+2\gamma^2}}\right)\right)}.$$

**5.3. Dual Smarandache $\tilde{t}\tilde{g}$ -curve of a unit dual spherical curve(ruled surface)**

In this section, we define the second type of dual Smarandache curves as dual Smarandache $\tilde{t}\tilde{g}$ -curve. Then, we give the relationships between the dual curve and its dual smarandache $\tilde{t}\tilde{g}$ -curve. Using the found results and relationships we study the developable of the corresponding ruled surface and its Smaranadache ruled surface.

**Definition 5.3.** Let $\tilde{\alpha}(\bar{s}) = \tilde{\alpha}$ be a unit speed regular dual curve lying fully on dual unit sphere and $\{\tilde{e}, \tilde{t}, \tilde{g}\}$ be its moving Darboux frame. The dual curve $\tilde{\alpha}_3$ defined by

$$\tilde{\alpha}_3 = \frac{1}{\sqrt{2}}\left(\tilde{t} + \tilde{g}\right)$$

is called the dual Smarandache $\tilde{t}\tilde{g}$ -curve of $\tilde{\alpha}$ and fully lies on $\tilde{S}^2$. Then the ruled surface corresponding to $\tilde{\alpha}_3$ is called the Smarandache $\vec{t}\vec{g}$ -ruled surface of the surface corresponding to dual curve $\tilde{\alpha}$.

Now we can give the relationships between $\tilde{\alpha}$ and its dual Smarandache $\tilde{t}\tilde{g}$ -curve $\tilde{\alpha}_3$ as follows.

***Theorem 5.13.*** *Let $\tilde{\alpha}(\bar{s}) = \tilde{\alpha}$ be a unit speed regular dual curve lying on dual unit sphere $\tilde{S}^2$. Then the relationships between the dual Darboux frames of $\tilde{\alpha}$ and its dual Smarandache $\tilde{t}\tilde{g}$ -curve $\tilde{\alpha}_3$ are given by*

$$\begin{pmatrix} \tilde{e}_3 \\ \tilde{t}_3 \\ \tilde{g}_3 \end{pmatrix} = \begin{pmatrix} 0 & \frac{1}{\sqrt{2}} & \frac{1}{\sqrt{2}} \\ \frac{-1}{\sqrt{1+2\bar{\gamma}^2}} & \frac{-\bar{\gamma}}{\sqrt{1+2\bar{\gamma}^2}} & \frac{\bar{\gamma}}{\sqrt{1+2\bar{\gamma}^2}} \\ \frac{2\bar{\gamma}}{\sqrt{2+4\bar{\gamma}^2}} & \frac{-1}{\sqrt{2+4\bar{\gamma}^2}} & \frac{1}{\sqrt{2+4\bar{\gamma}^2}} \end{pmatrix} \begin{pmatrix} \tilde{e} \\ \tilde{t} \\ \tilde{g} \end{pmatrix}. \tag{37}$$

If we represent the dual darboux frames of $\tilde{\alpha}$ and $\tilde{\alpha}_3$ by the dual matrixes $\tilde{E}$ and $\tilde{E}_1$, respectively, then (37) can be written as follows

$$\tilde{E} = \tilde{A}\tilde{E}_1$$

where

$$\tilde{A}=\begin{pmatrix} 0 & \frac{1}{\sqrt{2}} & \frac{1}{\sqrt{2}} \\ \frac{-1}{\sqrt{1+2\overline{\gamma}^2}} & \frac{-\overline{\gamma}}{\sqrt{1+2\overline{\gamma}^2}} & \frac{\overline{\gamma}}{\sqrt{1+2\overline{\gamma}^2}} \\ \frac{2\overline{\gamma}}{\sqrt{2+4\overline{\gamma}^2}} & \frac{-1}{\sqrt{2+4\overline{\gamma}^2}} & \frac{1}{\sqrt{2+4\overline{\gamma}^2}} \end{pmatrix}. \tag{38}$$

It is easily seen that $\det(\tilde{A})=1$ and $\tilde{A}\tilde{A}^T=\tilde{A}^T\tilde{A}=I$ where $I$ is the $3\times 3$ unitary matrix. It means that $\tilde{A}$ is a dual orthogonal matrix. Then we can give the following corollary.

***Corollary 5.7.*** *The relationship between the Darboux frames of the dual curves(ruled surfaces) $\tilde{\alpha}$ and $\tilde{\alpha}_3$ is given by a dual orthogonal matrix defined by (38).*

***Theorem 5.14.*** *The relationship between the dual Darboux formulae of dual Smarandache $\tilde{t}\tilde{g}$ -curve $\tilde{\alpha}_3$ and dual Darboux frame of $\tilde{\alpha}$ is given by*

$$\begin{pmatrix} \frac{d\tilde{e}_3}{ds_3} \\ \frac{d\tilde{t}_3}{ds_3} \\ \frac{d\tilde{g}_3}{ds_3} \end{pmatrix}=\begin{pmatrix} \frac{-1}{\sqrt{1+2\overline{\gamma}^2}} & \frac{-\overline{\gamma}}{\sqrt{1+2\overline{\gamma}^2}} & \frac{\overline{\gamma}}{\sqrt{1+2\overline{\gamma}^2}} \\ \frac{2\sqrt{2}\overline{\gamma}\overline{\gamma}'-\sqrt{2}\overline{\gamma}\left(1+2\overline{\gamma}^2\right)}{\left(1+2\overline{\gamma}^2\right)^2} & \frac{2\sqrt{2}\overline{\gamma}^2\overline{\gamma}'-\sqrt{2}\left(\overline{\gamma}'+\overline{\gamma}^2+1\right)\left(1+2\overline{\gamma}^2\right)}{\left(1+2\overline{\gamma}^2\right)^2} & \frac{-2\sqrt{2}\overline{\gamma}^2\overline{\gamma}'+\sqrt{2}\left(\overline{\gamma}'-\overline{\gamma}^2\right)\left(1+2\overline{\gamma}^2\right)}{\left(1+2\overline{\gamma}^2\right)^2} \\ \frac{-16\overline{\gamma}^2\overline{\gamma}'+(4\overline{\gamma}'+2)\left(2+4\overline{\gamma}^2\right)}{\left(1+2\overline{\gamma}^2\right)^2} & \frac{8\overline{\gamma}\overline{\gamma}'+2\overline{\gamma}\left(2+4\overline{\gamma}^2\right)}{\left(1+2\overline{\gamma}^2\right)^2} & \frac{-8\overline{\gamma}\overline{\gamma}'-2\overline{\gamma}\left(2+4\overline{\gamma}^2\right)}{\left(1+2\overline{\gamma}^2\right)^2} \end{pmatrix}\begin{pmatrix} \tilde{e} \\ \tilde{t} \\ \tilde{g} \end{pmatrix}.$$

(39)

***Theorem 5.15.*** *Let $\tilde{\alpha}(\overline{s})=\tilde{\alpha}$ be a unit speed regular curve on dual unit sphere. Then the relationship between the dual curvatures of $\tilde{\alpha}$ and its dual Smarandache $\tilde{t}\tilde{g}$ -curve $\tilde{\alpha}_3$ is given by*

$$\overline{\gamma}_3=\frac{4\sqrt{2}\overline{\gamma}\overline{\gamma}'+4\sqrt{2}\overline{\gamma}^2+2\sqrt{2}}{\left(2+4\overline{\gamma}^2\right)^{\frac{3}{2}}}.$$

***Corollary 5.8.*** *If the dual curvature $\overline{\gamma}$ of $\tilde{\alpha}$ is zero, the dual curvature $\overline{\gamma}_3$ of dual Smarandache $\tilde{e}\tilde{t}$ -curve $\tilde{\alpha}_3$ is* 1.

***Corollary 5.9.*** *The Darboux instantaneous vector of dual Smarandache $\tilde{t}\tilde{g}$ -curve is given by*

$$\tilde{d}_3=\frac{2\overline{\gamma}}{\sqrt{2+4\overline{\gamma}^2}}\,\tilde{e}+\frac{4\overline{\gamma}\overline{\gamma}'}{\left(2+4\overline{\gamma}^2\right)^{\frac{3}{2}}}\tilde{t}+\frac{4\overline{\gamma}\overline{\gamma}'+8\overline{\gamma}^2+4}{\left(2+4\overline{\gamma}^2\right)^{\frac{3}{2}}}\,\tilde{g}.$$

***Theorem 5.16.*** *Let $\tilde{\alpha}(\overline{s})=\tilde{\alpha}$ be a unit speed regular curve on dual unit sphere and $\tilde{\alpha}_3$ be the dual Smarandache $\tilde{t}\tilde{g}$ -curve of $\tilde{\alpha}$. If the ruled surface corresponding to the dual curve $\tilde{\alpha}$ is developable then the ruled surface corresponding to dual curve $\tilde{\alpha}_3$ is developable if and only if*

$$\delta_1 = \frac{4\sqrt{2}\,\delta\gamma' + 4\sqrt{2}\gamma\delta' + 8\sqrt{2}\delta\gamma}{\left(2+\gamma^2\right)^{\frac{3}{2}}} + \frac{12\delta\gamma(4\sqrt{2}\gamma\gamma' + 4\sqrt{2}\gamma^2 + 2\sqrt{2})}{\left(2+\gamma^2\right)^{\frac{5}{2}}}.$$

***Theorem 5.17.*** *Let $\tilde{\alpha}(\overline{s}) = \tilde{\alpha}$ be a unit speed regular curve on dual unit sphere. Then the relationship between the radius of dual curvature of dual Smarandache $\tilde{t}\tilde{g}$ -curve $\tilde{\alpha}_3$ and the dual curvature of $\tilde{\alpha}(\overline{s}) = \tilde{\alpha}$ is,*

$$\overline{R}_3 = \frac{\left(2+4\overline{\gamma}^2\right)^{\frac{3}{2}}}{\sqrt{(2+4\overline{\gamma}^2)^3 + (4\sqrt{2}\overline{\gamma}\overline{\gamma}' + 4\sqrt{2}\overline{\gamma}^2 + 2\sqrt{2})^2}}.$$

***Theorem 5.18.*** *Let $\tilde{\alpha}(\overline{s}) = \tilde{\alpha}$ be a unit speed regular curve on dual unit sphere. Then the relationship between the radius of dual spherical curvature of dual Smarandache $\tilde{t}\tilde{g}$ -curve $\tilde{\alpha}_3$ and the elements of dual curvature of $\tilde{\alpha}$ is,*

$$\overline{\rho}_3 = \arcsin\left(\frac{\left(2+4\overline{\gamma}^2\right)^{\frac{3}{2}}}{\sqrt{(2+4\overline{\gamma}^2)^3 + (4\sqrt{2}\overline{\gamma}\overline{\gamma}' + 4\sqrt{2}\overline{\gamma}^2 + 2\sqrt{2})^2}}\right).$$

### 5.4. Dual Smarandache $\tilde{e}\tilde{t}\tilde{g}$ -curve of a unit dual spherical curve(ruled surface)

In this section, we define the second type of dual Smarandache curves as dual Smarandache $\tilde{e}\tilde{t}\tilde{g}$ -curve. Then, we give the relationships between the dual curve and its dual smarandache $\tilde{e}\tilde{t}\tilde{g}$ -curve. Using the found results and relationships we study the developable of the corresponding ruled surface and its Smaranadache ruled surface.

**Definition 5.4.** Let $\tilde{\alpha}(\overline{s}) = \tilde{\alpha}$ be a unit speed regular dual curve lying fully on dual unit sphere and $\{\tilde{e}, \tilde{t}, \tilde{g}\}$ be its moving Darboux frame. The dual curve $\tilde{\alpha}_4$ defined by

$$\tilde{\alpha}_4 = \frac{1}{\sqrt{3}}\left(\tilde{e} + \tilde{t} + \tilde{g}\right)$$

is called the dual Smarandache $\tilde{e}\tilde{t}\tilde{g}$ -curve of $\tilde{\alpha}$ and fully lies on $\tilde{S}^2$. Then the ruled surface corresponding to $\tilde{\alpha}_4$ is called the Smarandache $\vec{e}\vec{t}\vec{g}$ -ruled surface of the surface corresponding to dual curve $\tilde{\alpha}$.

Now we can give the relationships between $\tilde{\alpha}$ and its dual smarandache $\tilde{e}\tilde{t}\tilde{g}$ -curve $\tilde{\alpha}_4$ as follows.

***Theorem 5.19.*** *Let $\tilde{\alpha}(\overline{s}) = \tilde{\alpha}$ be a unit speed regular dual curve lying on dual unit sphere $\tilde{S}^2$. Then the relationships between the dual Darboux frames of $\tilde{\alpha}$ and its dual Smarandache $\tilde{e}\tilde{t}\tilde{g}$ -curve $\tilde{\alpha}_4$ are given by*

$$\begin{pmatrix} \tilde{e}_4 \\ \tilde{t}_4 \\ \tilde{g}_4 \end{pmatrix} = \begin{pmatrix} \frac{1}{\sqrt{3}} & \frac{1}{\sqrt{3}} & \frac{1}{\sqrt{3}} \\ \frac{-1}{\sqrt{2\bar{\gamma}^2 - 2\bar{\gamma} + 2}} & \frac{1-\tilde{\gamma}}{\sqrt{2\bar{\gamma}^2 - 2\bar{\gamma} + 2}} & \frac{\tilde{\gamma}}{\sqrt{2\bar{\gamma}^2 - 2\bar{\gamma} + 2}} \\ \frac{2\bar{\gamma} - 1}{\sqrt{6}\sqrt{\bar{\gamma}^2 - \bar{\gamma} + 1}} & \frac{-(\bar{\gamma}+1)}{\sqrt{6}\sqrt{\bar{\gamma}^2 - \bar{\gamma} + 1}} & \frac{2-\bar{\gamma}}{\sqrt{6}\sqrt{\bar{\gamma}^2 - \bar{\gamma} + 1}} \end{pmatrix} \begin{pmatrix} \tilde{e} \\ \tilde{t} \\ \tilde{g} \end{pmatrix}. \tag{40}$$

If we represent the dual darboux frames of $\tilde{\alpha}$ and $\tilde{\alpha}_4$ by the dual matrixes $\tilde{E}$ and $\tilde{E}_1$, respectively, then (40) can be written as follows

$$\tilde{E} = \tilde{A}\tilde{E}_1$$

where

$$\tilde{A} = \begin{pmatrix} \frac{1}{\sqrt{3}} & \frac{1}{\sqrt{3}} & \frac{1}{\sqrt{3}} \\ \frac{-1}{\sqrt{2\bar{\gamma}^2 - 2\bar{\gamma} + 2}} & \frac{1-\tilde{\gamma}}{\sqrt{2\bar{\gamma}^2 - 2\bar{\gamma} + 2}} & \frac{\tilde{\gamma}}{\sqrt{2\bar{\gamma}^2 - 2\bar{\gamma} + 2}} \\ \frac{2\bar{\gamma} - 1}{\sqrt{6}\sqrt{\bar{\gamma}^2 - \bar{\gamma} + 1}} & \frac{-(\bar{\gamma}+1)}{\sqrt{6}\sqrt{\bar{\gamma}^2 - \bar{\gamma} + 1}} & \frac{2-\bar{\gamma}}{\sqrt{6}\sqrt{\bar{\gamma}^2 - \bar{\gamma} + 1}} \end{pmatrix}. \tag{41}$$

It is easily seen that $\det(\tilde{A}) = 1$ and $\tilde{A}\tilde{A}^T = \tilde{A}^T\tilde{A} = I$ where $I$ is the $3\times 3$ unitary matrix. It means that $\tilde{A}$ is a dual orthogonal matrix. Then we can give the following corollary.

***Corollary 5.10.*** *The relationship between the Darboux frames of the dual curves(ruled surfaces) $\tilde{\alpha}$ and $\tilde{\alpha}_4$ is given by a dual orthogonal matrix defined by (41).*

***Theorem 5.20.*** *The relationship between the dual Darboux formulae of dual Smarandache $\tilde{e}\tilde{t}\tilde{g}$-curve $\tilde{\alpha}_4$ and dual Darboux frame of $\tilde{\alpha}$ is given by*

$$\begin{pmatrix} \frac{d\tilde{e}_4}{ds_4} \\ \frac{d\tilde{t}_4}{ds_4} \\ \frac{d\tilde{g}_4}{ds_4} \end{pmatrix} = \begin{pmatrix} \frac{-1}{\sqrt{2\bar{\gamma}^2 - 2\bar{\gamma}+2}} & \frac{1-\bar{\gamma}}{\sqrt{2\bar{\gamma}^2 - 2\bar{\gamma}+2}} & \frac{\bar{\gamma}}{\sqrt{2\bar{\gamma}^2 - 2\bar{\gamma}+2}} \\ \frac{\sqrt{3}\bar{\gamma}'(2\bar{\gamma}-1)+\sqrt{3}(\bar{\gamma}-1)\left(2\bar{\gamma}^2-2\bar{\gamma}+2\right)^2}{\left(2\bar{\gamma}^2-2\bar{\gamma}+2\right)^3} & \frac{-\sqrt{3}\left(\bar{\gamma}'+\bar{\gamma}^2\right)\left(2\bar{\gamma}^2-2\bar{\gamma}+2\right)^2-\sqrt{3}(\bar{\gamma}-1)(2\bar{\gamma}\bar{\gamma}'-\bar{\gamma}')}{\left(2\bar{\gamma}^2-2\bar{\gamma}+2\right)^3} & \frac{\sqrt{3}\left(\bar{\gamma}'+\bar{\gamma}-\bar{\gamma}^2\right)\left(2\bar{\gamma}^2-2\bar{\gamma}+2\right)^2-\sqrt{3}\bar{\gamma}'\left(2\bar{\gamma}^2-\bar{\gamma}\right)}{\left(2\bar{\gamma}^2-2\bar{\gamma}+2\right)^3} \\ \frac{(4\bar{\gamma}'+2\bar{\gamma}+2)\left(\bar{\gamma}^2-\bar{\gamma}+1\right)-(2\bar{\gamma}-1)(2\bar{\gamma}\bar{\gamma}'-\bar{\gamma}')}{4\left(\bar{\gamma}^2-\bar{\gamma}+1\right)^2} & \frac{\left(-2\bar{\gamma}'+2\bar{\gamma}^2-2\right)\left(\bar{\gamma}^2-\bar{\gamma}+1\right)+(\bar{\gamma}+1)(2\bar{\gamma}\bar{\gamma}'-\bar{\gamma}')}{4\left(\bar{\gamma}^2-\bar{\gamma}+1\right)^2} & \frac{\left(-2\bar{\gamma}'-2\bar{\gamma}^2-2\right)\left(\bar{\gamma}^2-\bar{\gamma}+1\right)+(\bar{\gamma}-2)(2\bar{\gamma}\bar{\gamma}'-\bar{\gamma}')}{4\left(\bar{\gamma}^2-\bar{\gamma}+1\right)^2} \end{pmatrix} \begin{pmatrix} \tilde{e} \\ \tilde{t} \\ \tilde{g} \end{pmatrix}.$$

(42)

***Theorem 5.21.*** *Let $\tilde{\alpha}(\bar{s}) = \tilde{\alpha}$ be a unit speed regular curve on dual unit sphere. Then the relationship between the dual curvatures of $\tilde{\alpha}$ and its dual Smarandache $\tilde{e}\tilde{t}\tilde{g}$-curve $\tilde{\alpha}_4$ is given by*

$$\bar{\gamma}_4 = \frac{3\bar{\gamma}' + 2\bar{\gamma}^3 + 2}{2\sqrt{2}\left(\bar{\gamma}^2 - \bar{\gamma} + 1\right)^{\frac{3}{2}}}.$$

***Corollary 5.11.*** *If the dual curvature $\bar{\gamma}$ of $\tilde{\alpha}$ is zero, then the dual curvature $\bar{\gamma}_4$ of dual Smarandache $\tilde{e}\tilde{t}\tilde{g}$-curve $\tilde{\alpha}_4$ is* $\dfrac{1}{\sqrt{2}}$.

***Corollary 5.12.*** *The Darboux instantaneous vector of dual Smarandache $\tilde{e}\tilde{t}\tilde{g}$-curve is given by*

$$\tilde{d}_4 = \frac{3\bar{\gamma}' + 6\bar{\gamma}^3 - 6\bar{\gamma}^2 + 6\bar{\gamma}}{2\sqrt{6}\left(\bar{\gamma}^2 - \bar{\gamma} + 1\right)^{\frac{3}{2}}}\tilde{e} + \frac{3\bar{\gamma}'}{2\sqrt{6}\left(\bar{\gamma}^2 - \bar{\gamma} + 1\right)^{\frac{3}{2}}}\tilde{t} + \frac{3\bar{\gamma}' + 6\bar{\gamma}^2 - 6\bar{\gamma} + 6}{2\sqrt{6}\left(\bar{\gamma}^2 - \bar{\gamma} + 1\right)^{\frac{3}{2}}}\tilde{g}$$

***Theorem 5.22.*** *Let $\tilde{\alpha}(\bar{s}) = \tilde{\alpha}$ be a unit speed regular curve on dual unit sphere and $\tilde{\alpha}_4$ be the dual Smarandache $\tilde{e}\tilde{t}\tilde{g}$-curve of $\tilde{\alpha}$. If the ruled surface corresponding to the dual curve $\tilde{\alpha}$ is developable then the ruled surface corresponding to dual curve $\tilde{\alpha}_4$ is developable if and only if*

$$\delta_4 = \frac{6\gamma^2\delta + 3\delta'}{2\sqrt{2}\left(\gamma^2 - \gamma + 1\right)^{\frac{3}{2}}} - \frac{3\delta\left(2\gamma - 1\right)\left(3\gamma' + 2\gamma^3 + 2\right)}{4\sqrt{2}\left(\gamma^2 - \gamma + 1\right)^{\frac{5}{2}}}.$$

***Theorem 5.23.*** *Let $\tilde{\alpha}(\bar{s}) = \tilde{\alpha}$ be a unit speed regular curve on dual unit sphere. Then the relationship between the radius of dual curvature of dual Smarandache $\tilde{e}\tilde{t}\tilde{g}$-curve $\tilde{\alpha}_4$ and the dual curvature of $\tilde{\alpha}(\bar{s}) = \tilde{\alpha}$ is,*

$$\bar{R}_4 = \frac{2\sqrt{2}\left(\bar{\gamma}^2 - \bar{\gamma} + 1\right)^{\frac{3}{2}}}{\sqrt{8\left(\bar{\gamma}^2 - \bar{\gamma} + 1\right)^3 + \left(3\bar{\gamma}' + 2\bar{\gamma}^3 + 2\right)^2}}.$$

***Theorem 5.24.*** *Let $\tilde{\alpha}(\bar{s}) = \tilde{\alpha}$ be a unit speed regular curve on dual unit sphere. Then the relationship between the radius of dual spherical curvature of dual Smarandache $\tilde{e}\tilde{t}\tilde{g}$-curve $\tilde{\alpha}_4$ and the elements of dual curvature of $\tilde{\alpha}$ is,*

$$\bar{\rho}_4 = \arcsin\left(\frac{2\sqrt{2}\left(\bar{\gamma}^2 - \bar{\gamma} + 1\right)^{\frac{3}{2}}}{\sqrt{8\left(\bar{\gamma}^2 - \bar{\gamma} + 1\right)^3 + \left(3\bar{\gamma}' + 2\bar{\gamma}^3 + 2\right)^2}}\right).$$

**Example 1.** Let consider the dual spherical curve $\tilde{\alpha}(\bar{s})$ given by the parametrization

$$\tilde{\alpha}(\bar{s}) = \left(\cos s, \sin s, 0\right) + \varepsilon\left(-s\sin s, s\cos s, 0\right).$$

The curve $\tilde{\alpha}(\bar{s})$ represents the ruled surface

$$r(s,v) = \left(v\cos s, v\sin s, s\right)$$

which is a helicoids surface rendered in Fig. 1. Then the dual Darboux frame of $\tilde{\alpha}$ is obtained as follows,

$$\tilde{e}(\bar{s}) = \left(\cos s, \sin s, 0\right) + \varepsilon\left(-s\sin s, s\cos s, 0\right)$$

$$\tilde{t}(\bar{s}) = \left(-\sin s, \cos s, 0\right) + \varepsilon\left(-s\cos s, -s\sin s, 0\right)$$

$$\tilde{g}(\bar{s}) = \left(0, 0, 1\right)$$

The Smarandache $\tilde{e}\tilde{t}$, $\tilde{e}\tilde{g}$, $\tilde{t}\tilde{g}$, and $\tilde{e}\tilde{t}\tilde{g}$ curves of the dual curve $\tilde{\alpha}$ are given by

$$\tilde{\alpha}_1(\bar{s}) = \frac{1}{\sqrt{2}}\left[(\cos s - \sin s, \cos s + \sin s, 0) + \varepsilon(-s\cos s - s\sin s, s\cos s - s\sin s, 0)\right]$$

$$\tilde{\alpha}_2(\bar{s}) = \frac{1}{\sqrt{2}}\left[(\cos s, \sin s, 1) + \varepsilon(-s\sin s, s\cos s, 0)\right]$$

$$\tilde{\alpha}_3(\bar{s}) = \frac{1}{\sqrt{2}}\left[(-\sin s, \cos s, 1) + \varepsilon(-s\cos s, -s\sin s, 0)\right]$$

$$\tilde{\alpha}_4(\bar{s}) = \frac{1}{\sqrt{3}}\left[(\cos s - \sin s, \cos s + \sin s, 1) + \varepsilon(-s\cos s - s\sin s, s\cos s - s\sin s, 0)\right]$$

respectively. From E. Study mapping, these dual spherical curves correspond to the following ruled surfaces

$$r_1(s,v) = (0,0,s) + v\left(\frac{1}{\sqrt{2}}\cos s - \frac{1}{\sqrt{2}}\sin s, \frac{1}{\sqrt{2}}\cos s + \frac{1}{\sqrt{2}}\sin s, 0\right)$$

$$r_2(s,v) = (0,0,s) + v\left(\frac{1}{\sqrt{2}}\cos s, \frac{1}{\sqrt{2}}\sin s, \frac{1}{\sqrt{2}}\right)$$

$$r_3(s,v) = (0,0,s) + v\left(-\frac{1}{\sqrt{2}}\sin s, \frac{1}{\sqrt{2}}\cos s, \frac{1}{\sqrt{2}}\right)$$

$$r_4(s,v) = (0,0,s) + v\left(\frac{1}{\sqrt{3}}\cos s - \frac{1}{\sqrt{3}}\sin s, \frac{1}{\sqrt{3}}\cos s + \frac{1}{\sqrt{3}}\sin s, \frac{1}{\sqrt{3}}\right)$$

respectively. These surfaces are rendered in Fig.2, Fig. 3, Fig. 4 and Fig. 5, respectively.

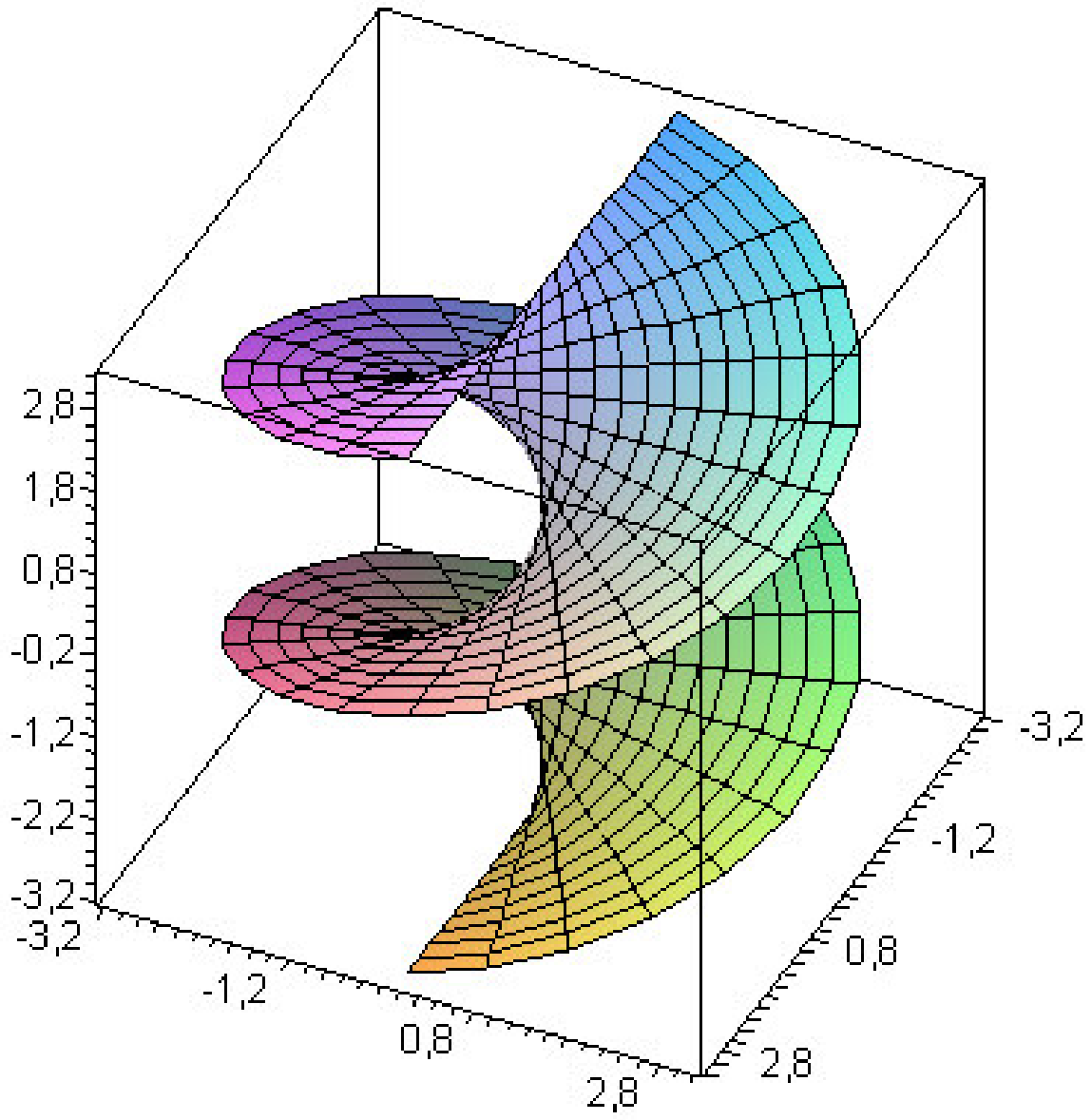

**Figure 1.** Helicoid surface corresponding to dual curve $\tilde{\alpha}$

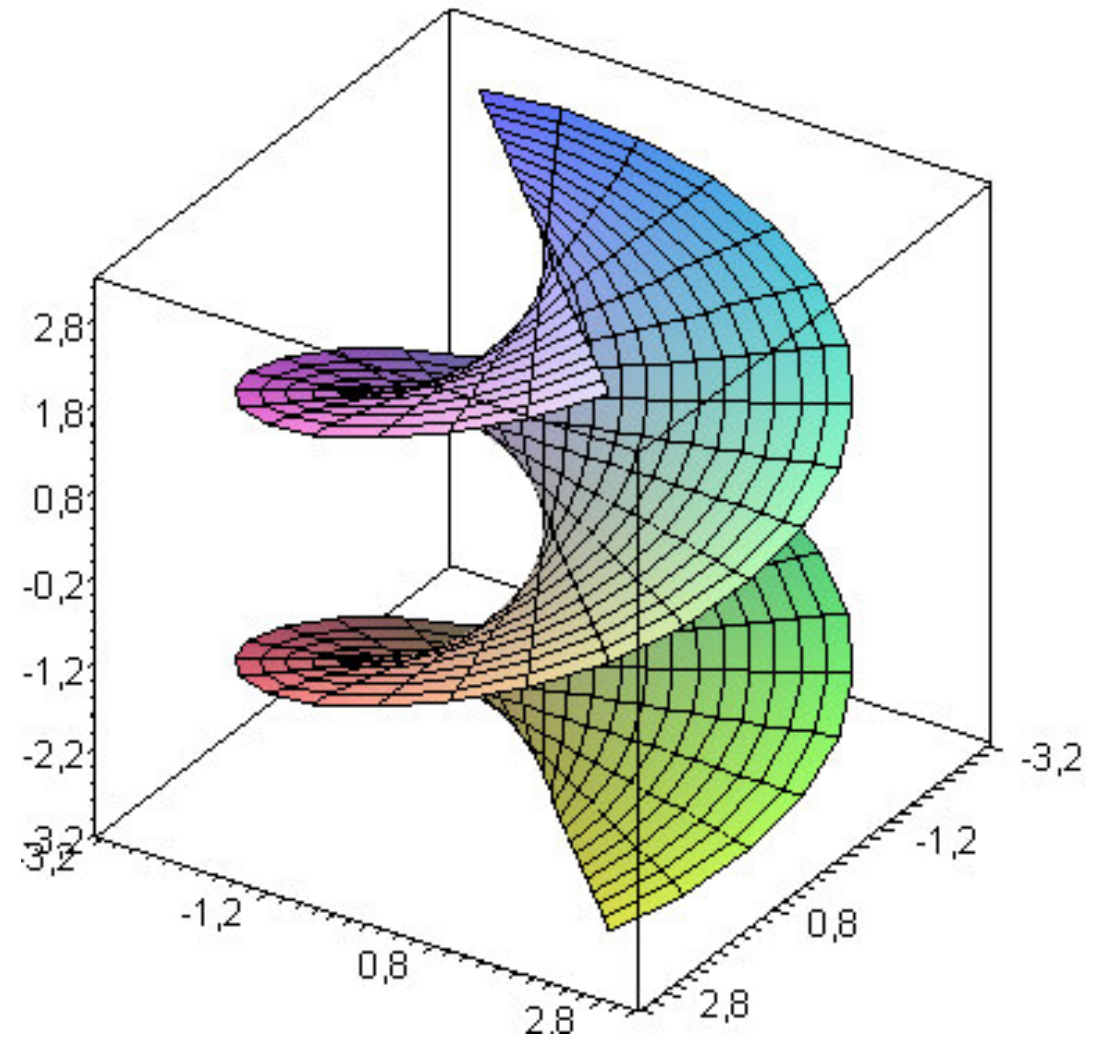

**Figure 2.** Smarandache $\tilde{e}\tilde{t}$ ruled surface

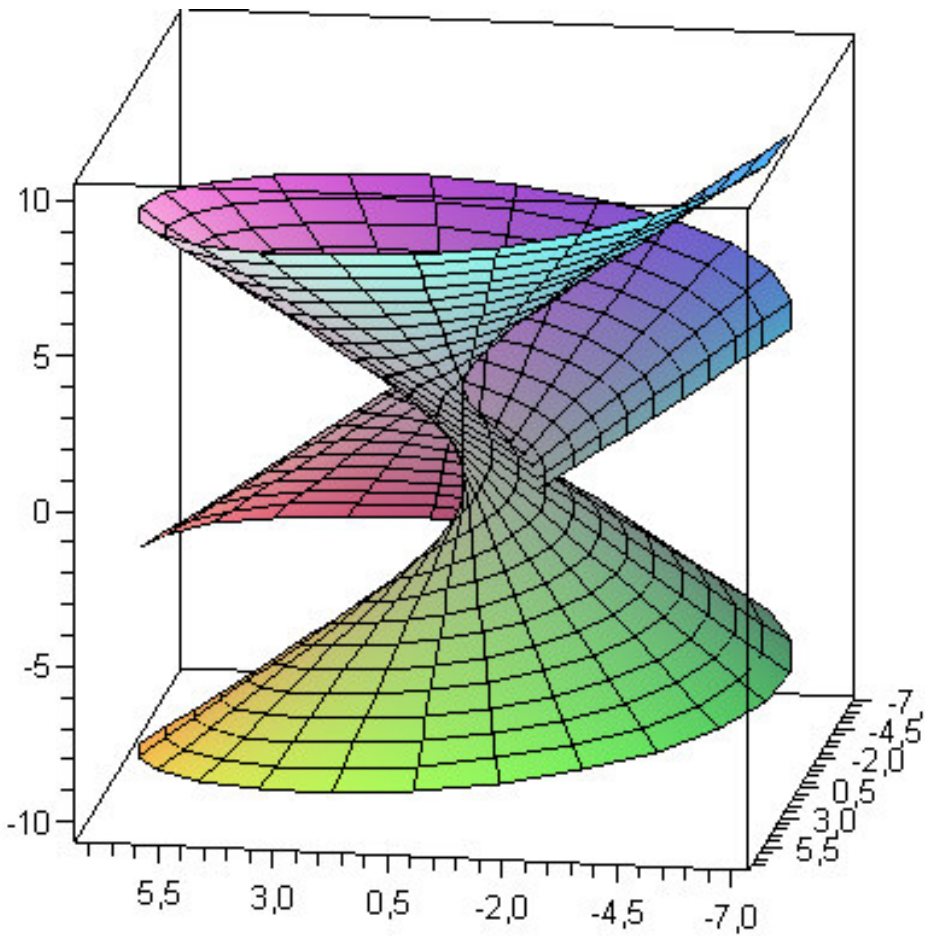

**Figure 3.** Smarandache $\tilde{e}\tilde{g}$ ruled surface

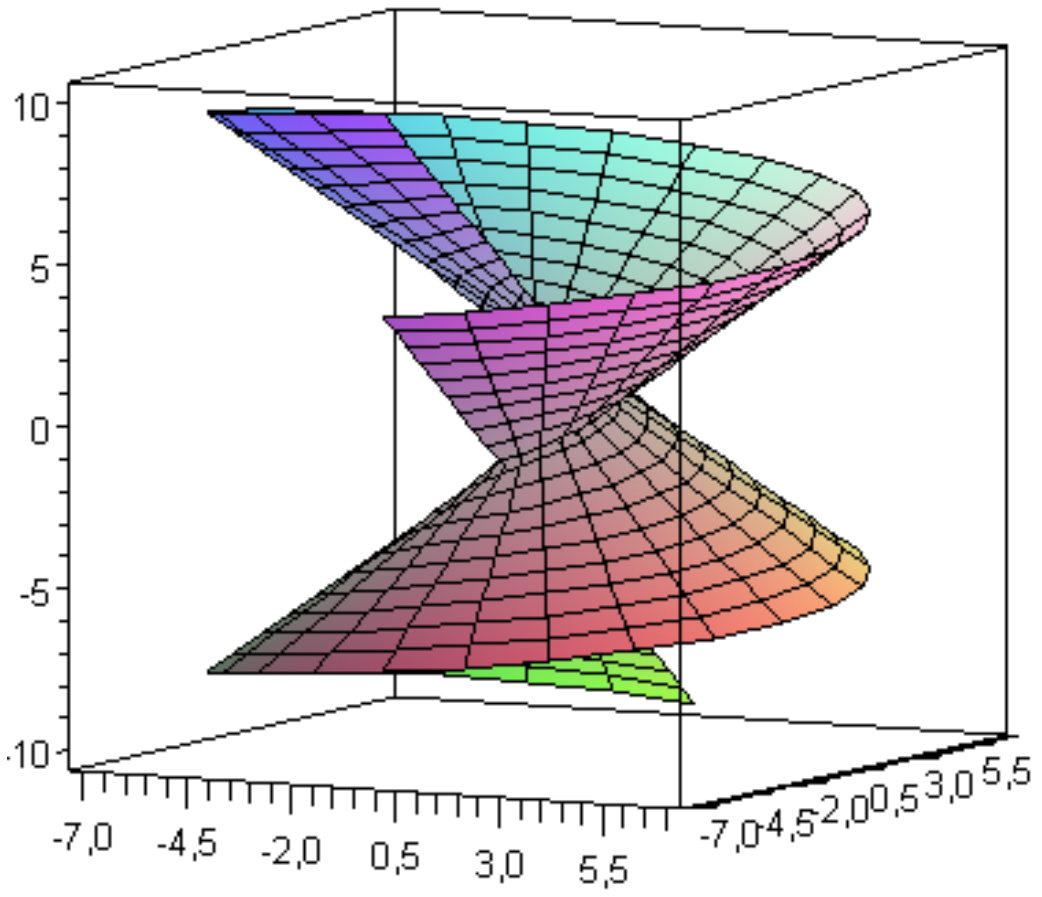

**Figure 4.** Smarandache $\tilde{t}\tilde{g}$ ruled surface

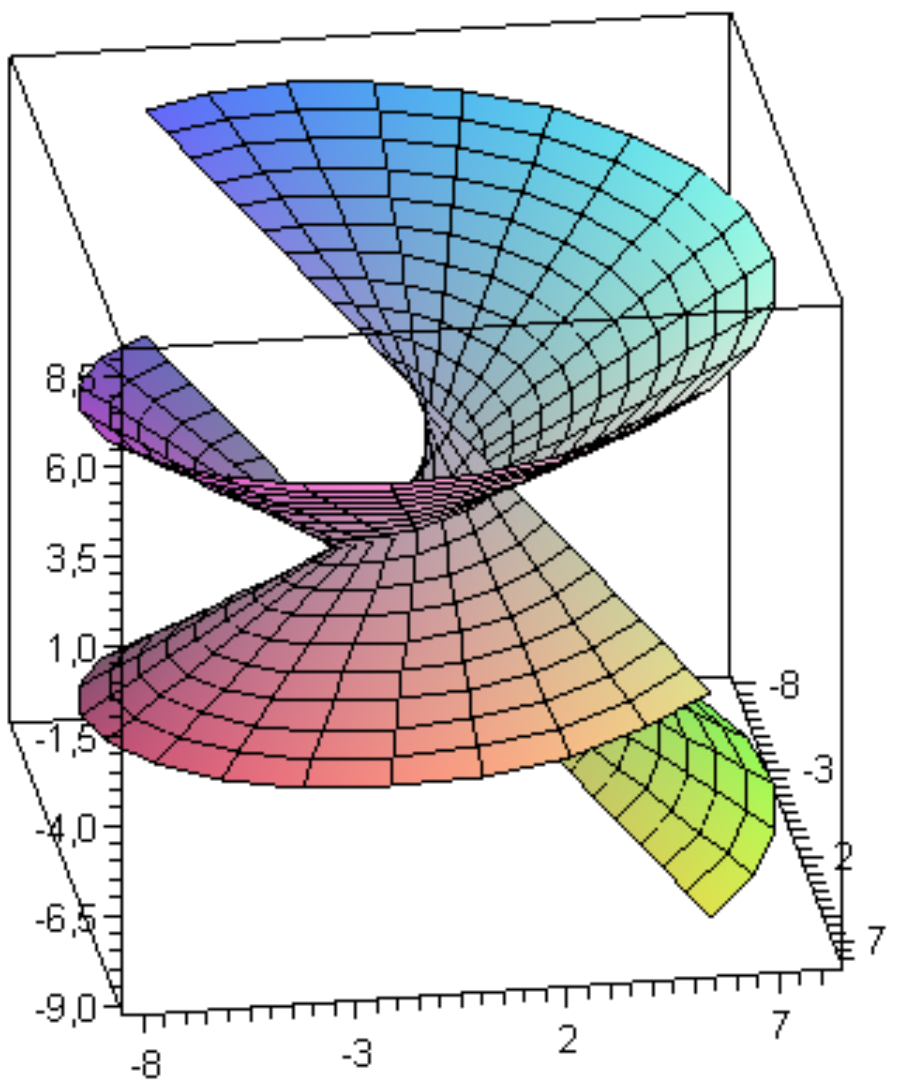

**Figure 5.** Smarandache $\tilde{e}\tilde{t}\tilde{g}$ ruled surface